# Asymptotic Expansions and Extrapolation of Approximate Eigenvalues for Second Order Elliptic Problems by Mixed Finite Element Methods

HEHU XIE*

**Abstract**

In this paper, we derive an asymptotic error expansion for the eigenvalue approximations by the lowest order Raviart-Thomas mixed finite element method for the general second order elliptic eigenvalue problems. Extrapolation based on such an expansion is applied to improve the accuracy of the eigenvalue approximations. Furthermore, we also prove the superclose property between the finite element projection with the finite element approximation of the eigenvalue problems by mixed finite element methods. In order to prove the full order of the eigenvalue extrapolation, we first propose " the auxiliary equation method". The result of this paper provides a general procedure to produce an asymptotic expansions for eigenvalue approximations by mixed finite elements.

**AMS subject classifications.** 80M10, 65L15, 65L60, 65L70, 65B99

*Academy of Mathematics and Systems Science, Chinese Academy of Sciences, Beijing 100080, China.